\documentclass[letterpaper]{article}
\usepackage{authblk}
\usepackage{url}
\usepackage[margin=1in,footskip=0.25in]{geometry}
\usepackage[labelfont=bf]{caption}
\usepackage{graphicx}
\usepackage{amsmath,amssymb,amsfonts,amsmath}
\usepackage{graphicx}
\usepackage{subfigure}
\usepackage{caption}
\usepackage{tikz}
\usepackage{float}
\usepackage[dvipdf]{epsfig}
\usepackage{float}
\usepackage{lipsum}
\usepackage{wrapfig}
\usepackage{cite}
\usepackage{color}
\usepackage{easylist}
\usepackage{lineno}
\usepackage{booktabs}
\usepackage{hyperref}
\usepackage{cleveref}
\usepackage{multirow}
\usepackage{hhline}
\usepackage{soul}
\usepackage{blindtext}
\usepackage{threeparttable}
\usepackage[bbgreekl]{mathbbol}

\usepackage{todonotes}

\def \L2{L^2}

\graphicspath{{Figures/}}

\newcommand{\eqn}[1]{
\begin{eqnarray}
#1
\end{eqnarray}
}

\newcommand{\tb}[1]{
{\bf #1}
}

\newcommand{\norm}[1]{
\lVert #1\rVert
}

\newcommand{\ip}[1]{
\langle #1\rangle
}

\newcommand{\bs}[1]{
\boldsymbol #1
}

\title{Improving the accuracy of physics-informed neural networks via last-layer retraining}

\author[1]{Saad Qadeer\thanks{{\tt saad.qadeer@pnnl.gov}}} 
\author[2]{Panos Stinis} 
\affil[1]{Advanced Computing, Mathematics and Data Division, Pacific Northwest National Laboratory, Richland, WA}
\affil[2]{Advanced Computing, Mathematics and Data Division, Pacific Northwest National Laboratory, Richland, WA}
\affil[ ]{Department of Applied Mathematics, University of Washington, Seattle, WA}

\date{\today}

\begin{document}

\maketitle

\begin{abstract}
Physics-informed neural networks (PINNs) are a versatile tool in the burgeoning field of scientific machine learning for solving partial differential equations (PDEs). However, determining suitable training strategies for them is not obvious, with the result that they typically yield moderately accurate solutions. In this article, we propose a method for improving the accuracy of PINNs by coupling them with a post-processing step that seeks the best approximation in a function space associated with the network. We find that our method yields errors four to five orders of magnitude lower than those of the parent PINNs across architectures and dimensions. Moreover, we can reuse the basis functions for the linear space in more complex settings, such as time-dependent and nonlinear problems, allowing for transfer learning. Our approach also provides a residual-based metric that allows us to optimally choose the number of basis functions employed.
\end{abstract}

\section{Introduction}\label{SecIntro}
The use of machine learning methods for problems in scientific computing has seen a huge surge in recent years. Driven chiefly by improvements in architectures, advances in training methods, and increasing ease of implementation, this interface (dubbed scientific machine learning) holds the potential of providing breakthroughs in long-standing challenges such as complex geometries, the curse of dimensionality, and integrating physics and experimental data. One of the most versatile tools pioneered in this field are physics-informed neural networks (PINNs), which rely on the universal approximation theorem to express solutions to partial differential equations (PDEs) as neural networks \cite{raissi2019physics}. Further, the use of automatic differentiation (AD) to train them to obey the governing systems imbues them with a simplicity that has led to widespread uptake across various domains.

At the same time, the study of limitations of PINNs has also drawn significant interest \cite{fuks2020limitations,wang2021understanding,mcgreivy2024weak}. The areas of inquiry range over improving the selection of training weights \cite{wang2022and,chen2024self} and placement of collocation points \cite{gao2023failure,tang2023pinns,chen2025self}, developing a better understanding of the loss landscapes \cite{krishnapriyan2021characterizing}, and modifying the architectures \cite{zeng2022competitive} and loss functions \cite{kharazmi2019variational,park2024beyond} to reflect the characteristics of the underlying problem. While these ideas have shown considerable promise in broadening the applicability of PINNs, they still await a comprehensive unified treatment to realize their full potential.

The flexibility of PINNs has made them particularly popular for providing low-fidelity approximations that can then be further refined by various techniques. Among others, these include multifidelity methods \cite{penwarden2022multifidelity,howard2023multifidelity}, multistage \cite{wang2024multi} and stacked \cite{howard2023stacked} architectures, and pre-trained finite element methods \cite{wang2026pretrain}. In a similar vein, PINNs can be sequentially trained to supply basis functions, the linear combinations of which can be used to compute solutions either in the Galerkin \cite{ainsworth2021galerkin,ainsworth2022galerkin} or least-squares sense \cite{weng2026deep}. Another highly promising avenue is seeking the solution in the linear space defined by the span of the last-layer neurons \cite{datar2024fast,jia2026physics}. The Physics-Driven Orthogonal Feature Method (PD-OFM) method, proposed in \cite{jia2026physics}, drives these neuron functions to be orthonormal to improve their expressivity, and uses them in an over-determined collocation formulation to solve the problem. The latter procedure is equivalent to retraining the last parameter layer and was proposed in \cite{qadeer2023efficient} as an inexpensive way of improving NN performance more broadly.

In this article, we propose using PINNs for identifying orthonormal bases for arbitrary spatial domains and using them in a spectral procedure to solve the given problem. Like the PD-OFM, our approach also employs the space defined by the neurons in the last hidden layer. However, rather than optimizing the same loss function as the trained PINN, it uses the basis functions in a variational formulation. The identification of spatial basis functions imparts our technique with greater flexibility, allowing them to be used seamlessly in more complex settings (echoing transfer learning) such as time-dependent and nonlinear problems on the same domain. Moreover, the basis extraction involves an explicit orthogonalization procedure that leads to smaller and better conditioned systems. The conditioning is also improved by our reliance on a variational formulation instead of a least-squares formulation, lowering the orders of derivatives involved. We find that our method yields errors four to five orders of magnitude lower than those of the parent PINNs, and also provides error control metrics by monitoring the residuals.

\section{Mathematical formulation for the Poisson equation}\label{SecApproach}

Let $\Omega \subset \mathbb{R}^d$ be open and connected domain. Consider the Poisson equation
\eqn{
-\nabla \cdot (k(x) \nabla u(x)) &=& f(x), \quad x \in \Omega, \nonumber\\
u(x) &=& g(x), \quad x \in \partial \Omega. \label{PoiProb}
}

Let $u_\text{NN}^{\theta}$ be a neural network trained to approximate the solution to \eqref{PoiProb}, using data and/or physics-informed losses. We propose a strategy that improves the accuracy of this solution by retraining the last layer of this architecture, as well as providing tools to solve more complex problems on $\Omega$.  

For concreteness, suppose $u_\text{NN}^{\theta}$ is a fully connected network (FCN) with $L$ hidden neuron layers of the form
\eqn{
\bs{z}^{(l)} &=& \sigma \left(W^{(l)}\bs{z}^{(l-1)} + \tb{b}^{(l)} \right), \quad 1 \leq l \leq L, \nonumber\\
\bs{z}^{(L+1)} &=& W^{(L+1)}\bs{z}^{(L)} + \tb{b}^{(L+1)}. \label{NN}
}

Here, $\sigma$ is the activation function, $W^{(l)} \in \mathbb{R}^{d_l \times d_{l-1}}$ and $\tb{b}^{(l)} \in \mathbb{R}^{d_l}$ are the trainable parameters, $\bs{z}^{(0)} = x \in \Omega$ is the input, and $\bs{z}^{(L+1)} \in \mathbb{R}$ is the output, i.e., $u_\text{NN}^{\theta}\left(\bs{z}^{(0)}\right) = \bs{z}^{(L+1)}$ (and $d_{L+1} = 1$). We denote all the trainable parameters $\{W^{(l)},\tb{b}^{(l)} \}_{1 \leq l \leq L+1}$ by $\theta$ for brevity.

The key feature of this architecture that we leverage is that no activation is applied at the last stage and hence the output is a linear combination of the last hidden layer of neurons. It follows that our approach is not limited only to this setup, and that any architecture that possesses this property is amenable to the following treatment.

Denote the neurons in the last hidden layer by $\{\phi_j\}_{1 \leq j \leq {d_L}}$ and set $\phi_0 \equiv 1$. Thus, \eqref{NN} can be rewritten as $u_\text{NN}^{\theta}(x) = \sum_{j = 0}^{d_L} w_j\phi_j(x)$, where $w_0 = \tb{b}^{(L+1)}$ and $w_j = W^{(L+1)}_j$. Note that since these coefficients have been found by stochastic gradient descent applied to all the trainable parameters (of which they are only a subset), they are possibly sub-optimal. Our approach is to improve the accuracy by finding the best solution to \eqref{PoiProb} in $\mathcal{S}_{\theta} = \text{span}\left(\{\phi_j\}_{j = 0}^{d_L}\right)$. 

Let $\{(x_i,\omega_i)\}_{1 \leq i \leq N_q}$ be a high-order quadrature rule on $\Omega$, i.e.,
\eqn{
\int_\Omega h(x) \ dx \approx \sum_{i = 1}^{N_q} h(x_i)\omega_i \nonumber
}
with high accuracy. Define $\Phi \in \mathbb{R}^{N_q \times (d_{L}+1)}$ by
\eqn{
\Phi_{ij} = \omega_i^{1/2}\phi_j(x_i), \label{PhiDefn}
}
and compute its singular-value decomposition (SVD) $\Phi = QSV^\top$. This provides a hierarchical orthonormal basis $\{q_k\}_{k = 0}^{d_L}$ for $\mathcal{S}_{\theta}$ that obeys $q_k(x_i) = \omega_i^{-1/2}Q_{ik}$ at the quadrature grid specifically, and
\eqn{
q_k(x) = \frac{1}{S_{kk}}\sum_{l = 0}^{d_L} \phi_l(x)V_{lk} \label{CoB}
}
more generally. Since the singular values decay rapidly, in practice we enforce a threshold $r \leq d_L$ so division by small singular values does not corrupt the corresponding basis functions and magnify small rounding errors (more on this later). We note that this procedure is similar to the one employed in \cite{meuris2023machine,williams2024physics} for extracting basis functions from a trained Deep Operator Network.

Next, in order to find the best approximation to \eqref{PoiProb} in $\mathcal{S}_{\theta}$, we use the Nitsche variational formulation as it allows us to seamlessly employ basis functions and arbitrary boundary data on complex domains \cite{nitsche1971variationsprinzip,benzaken2022constructing}. This approach solves \eqref{PoiProb} by posing the minimization problem $\min_{v \in \mathcal{S}} \mathcal{J}[v]$ over a space $\mathcal{S}$, where
\eqn{
\mathcal{J}[v] = \frac{1}{2}\int_{\Omega} k\left|\nabla v\right|^2 - fu \ dx - \int_{\partial \Omega} k \left(\partial_{\bf{n}}v \right) (v-g) \ ds + \beta \int_{\partial \Omega} k(v-g)^2 \ ds. \label{Nitsche}
}

Here, $\tb{n}$ is the outward facing normal vector, and $\beta > 0$ is chosen large enough to ensure that $\mathcal{J}$ is strictly convex and hence possesses a unique minimizer. In the results that follow, we set $\beta = 200$.

We can minimize this functional over $\mathcal{S}_{\theta}$ by writing $u_r(x) = \sum_{j = 0}^r c_jq_j(x)$, defining $J(\tb{c}) = \mathcal{J}[u_r]$, setting $\nabla_{\bf{c}} J = 0$ and solving the resulting linear system $A\tb{c} = \tb{b}$, where
\eqn{
A_{lj} &=& -\ip{\nabla q_l , k \nabla q_j} + \int_{\partial \Omega} k\left(q_l \left(\partial_{\bf{n}}q_j  \right)  + \left(\partial_{\bf{n}} q_l \right) q_j - \beta q_lq_j \right) \ ds, \nonumber\\
b_l &=& - \ip{q_l,f} + \int_{\partial \Omega} (\partial_{\bf{n}}q_l)g - \beta q_l g \ ds, \qquad 0 \leq l,j \leq r, \label{Abdefn}
}
and where $\ip{h_1,h_2} = \int_{\Omega} h_1(x)h_2(x) \ dx$. The derivatives of the basis functions can be calculated using \eqref{CoB} and making use of the structure of \eqref{NN} to compute the derivatives of $\{\phi_l\}$ exactly. We can also use it to calculate the residual $e_r(x) = -\nabla \cdot (k(x) \nabla u_r(x)) - f(x)$, the size of which can then be used to guide us in choosing the optimal value of $r$.

In order to test this approach, we consider \eqref{PoiProb} with $k(x) \equiv 1$ on 
\begin{itemize}
\item[(i)] an interval $\Omega_1 = (-1,1)$, with $f(x) = \sin(6x)$ and $g(-1) = 0$, $g(1) = 2$, and compare against the (easily calculated) analytical solution;

\item[(ii)] a square box $\Omega_2 = (-1,1)^2$, and use the exact solution $u(x,y) = \cos(\sin(x)-y^2)$ to supply the $f$ and $g$;

\item[(iii)] an L-shaped domain $\Omega_3 = (-1,1)^2 - [0,1)^2$ with the same exact solution as in (ii).

\end{itemize}

\begin{figure}[tbph]
\centering
\subfigure[$D = \{1,30,30,1\}$]
{\includegraphics[width=0.31\textwidth]{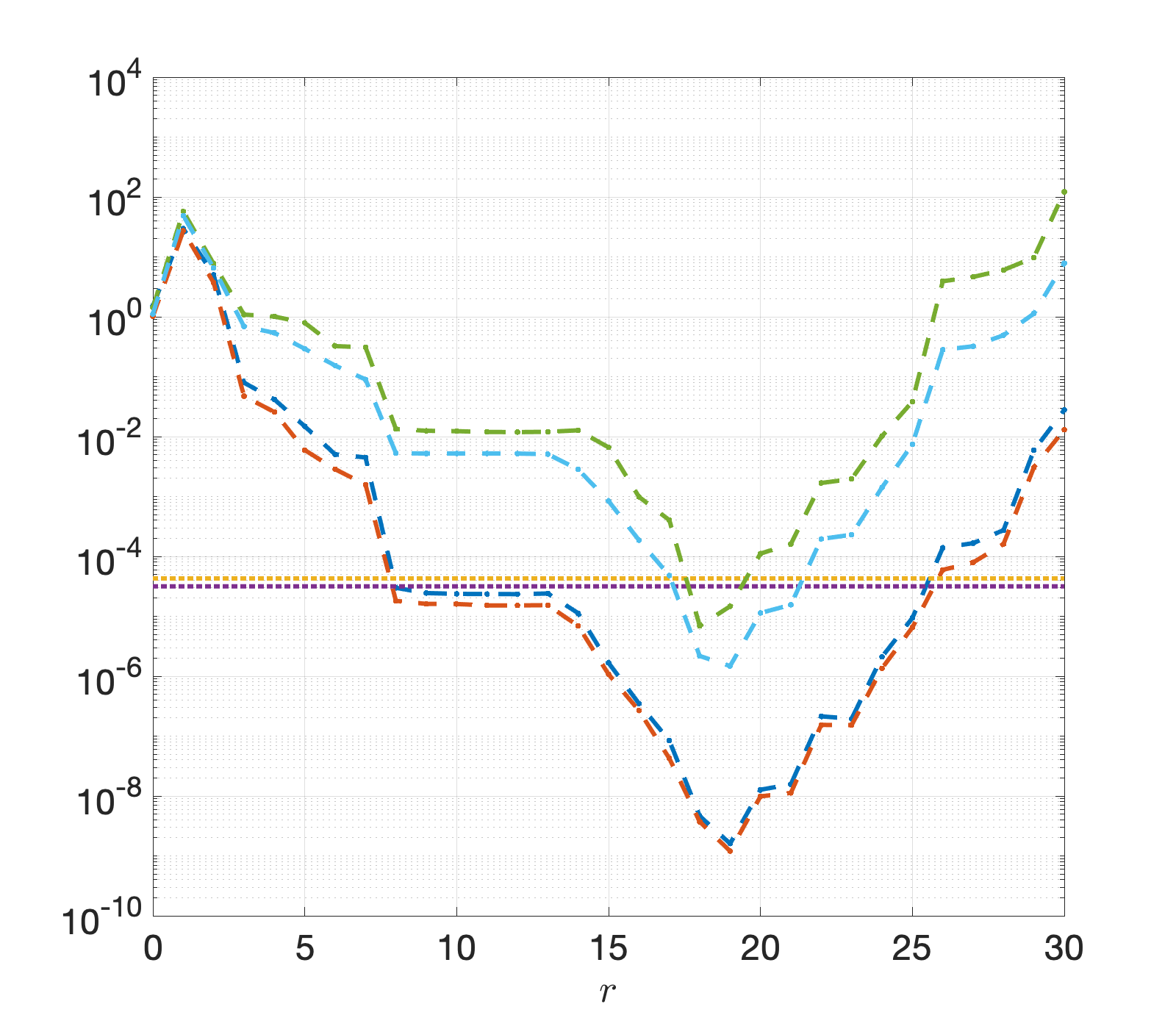}
}
\subfigure[$D = \{1,60,60,1\}$]
{\includegraphics[width=0.31\textwidth]{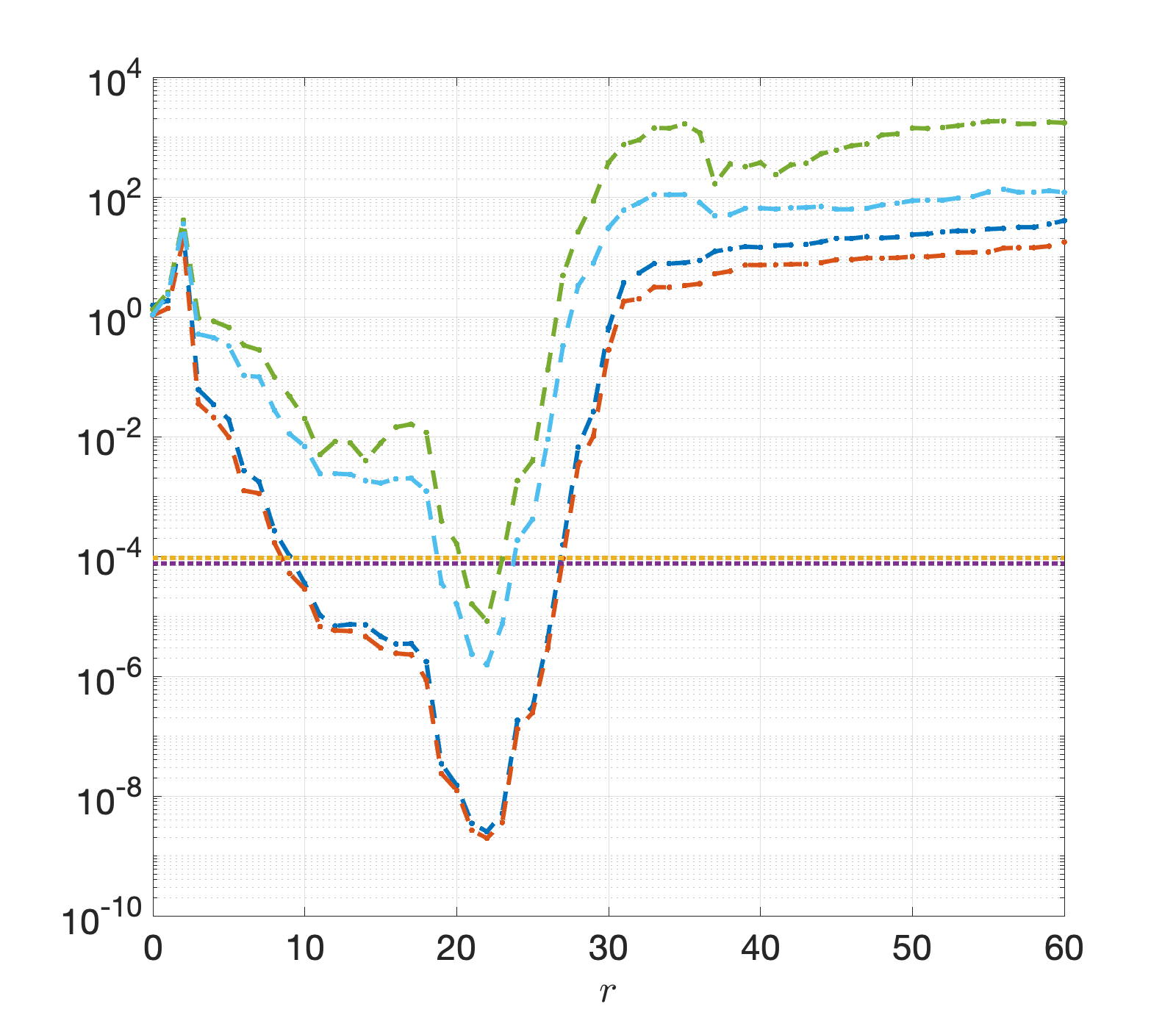}
}
\subfigure[$D = \{1,50,50,50,1\}$]
{\includegraphics[width=0.31\textwidth]{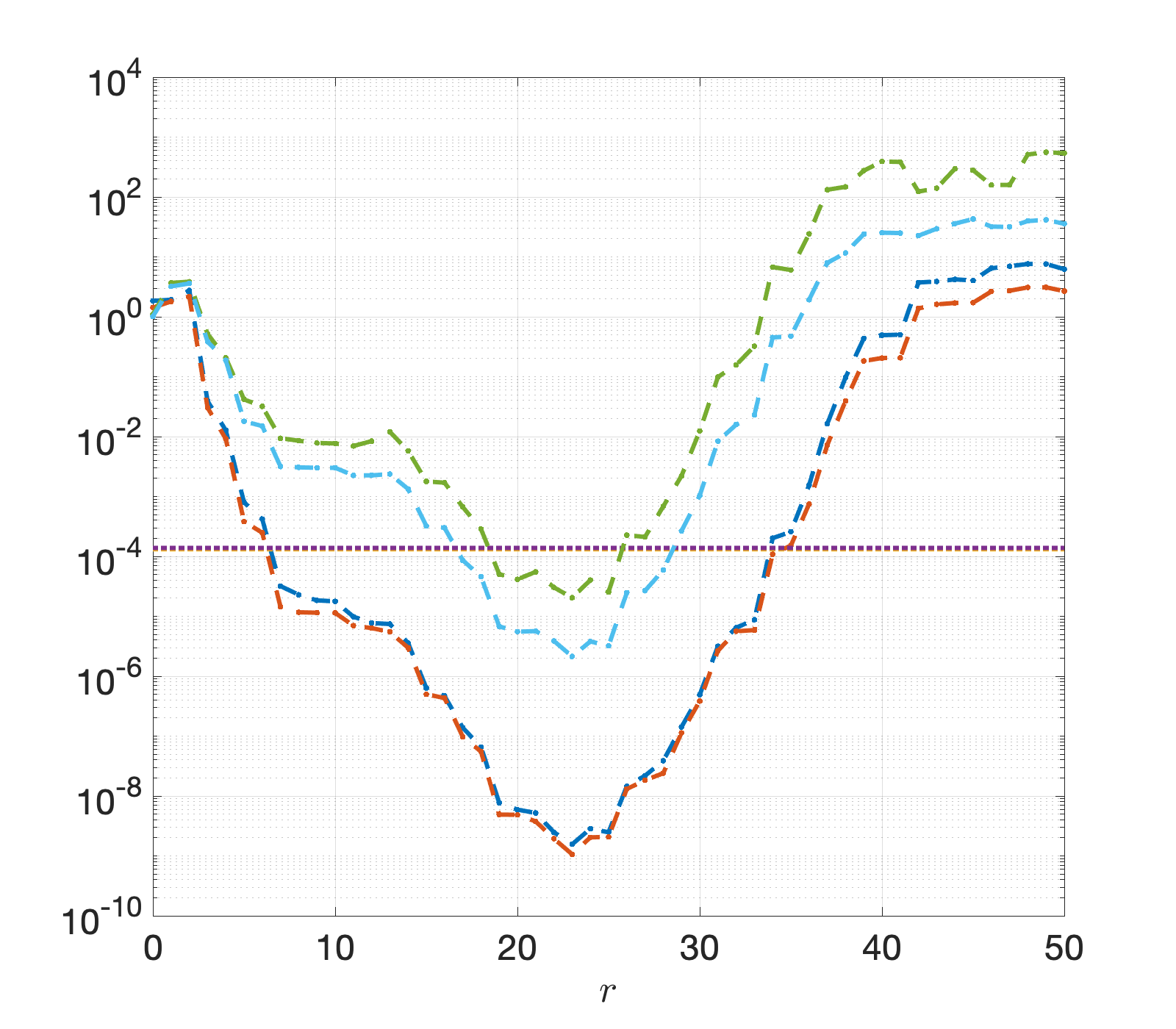}
}
\caption{Poisson equation on the one dimensional domain $\Omega_1$. The errors $\norm{u_r - u}$ for different values of $r$, measured in $L^{\infty}$ (blue) and $L^2$ (red), for problem (i), across PINN architectures. We also show the PINN errors $\norm{u^{\theta}_\text{NN} - u}$ in $L^{\infty}$ (yellow) and $L^2$ (purple) as the horizontal curves. The residuals $\norm{e_r}$ in in $L^{\infty}$ (green) and $L^2$ (cyan) are also shown; note the close similarity in the behaviour of the error and residual curves.}\label{PI_Poi1}
\end{figure}

\begin{figure}[tbph]
\centering
\subfigure[$D = \{2,200,400,1\}$]
{\includegraphics[width=0.31\textwidth]{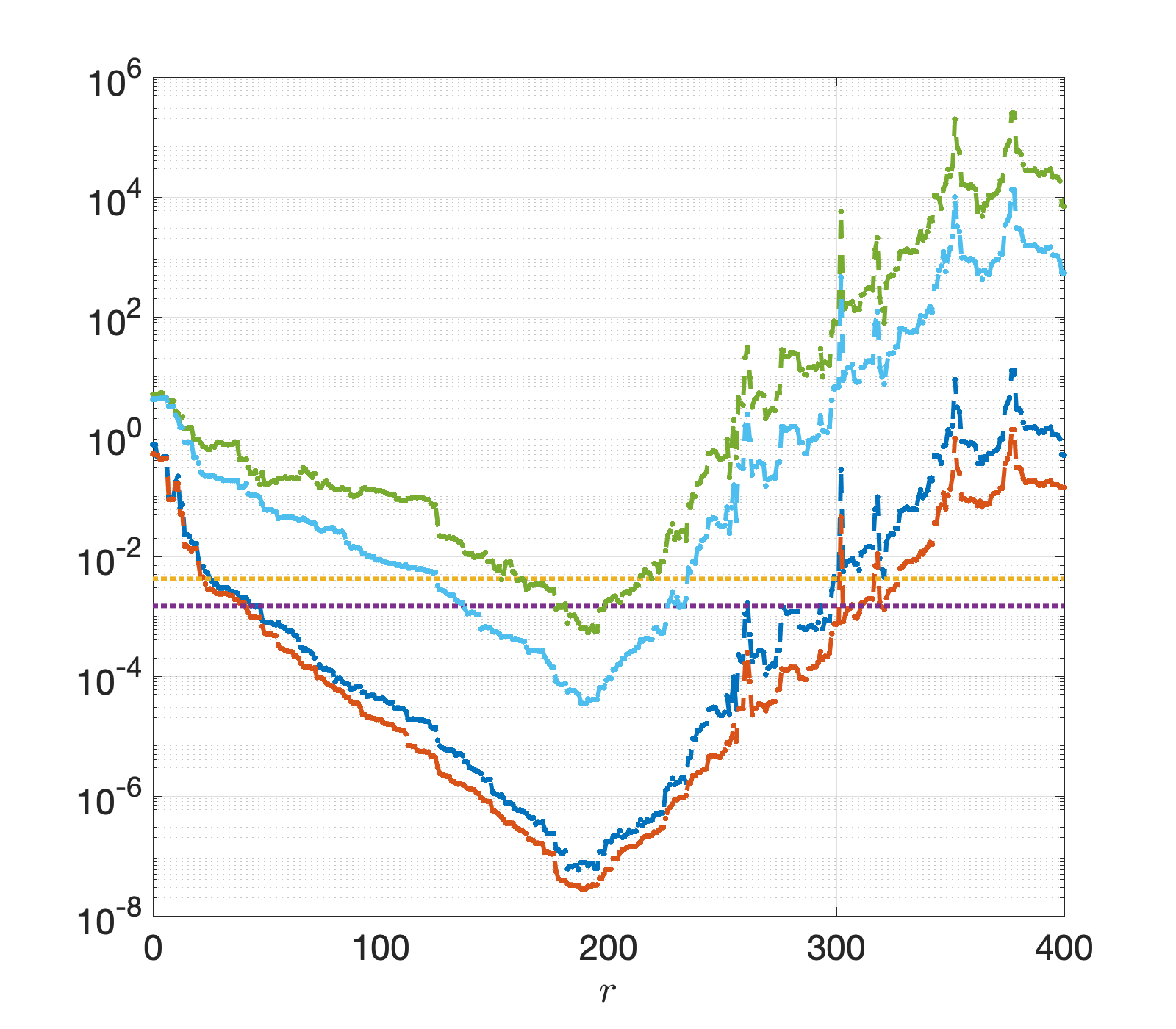}
}
\subfigure[$D = \{2,300,400,1\}$]
{\includegraphics[width=0.31\textwidth]{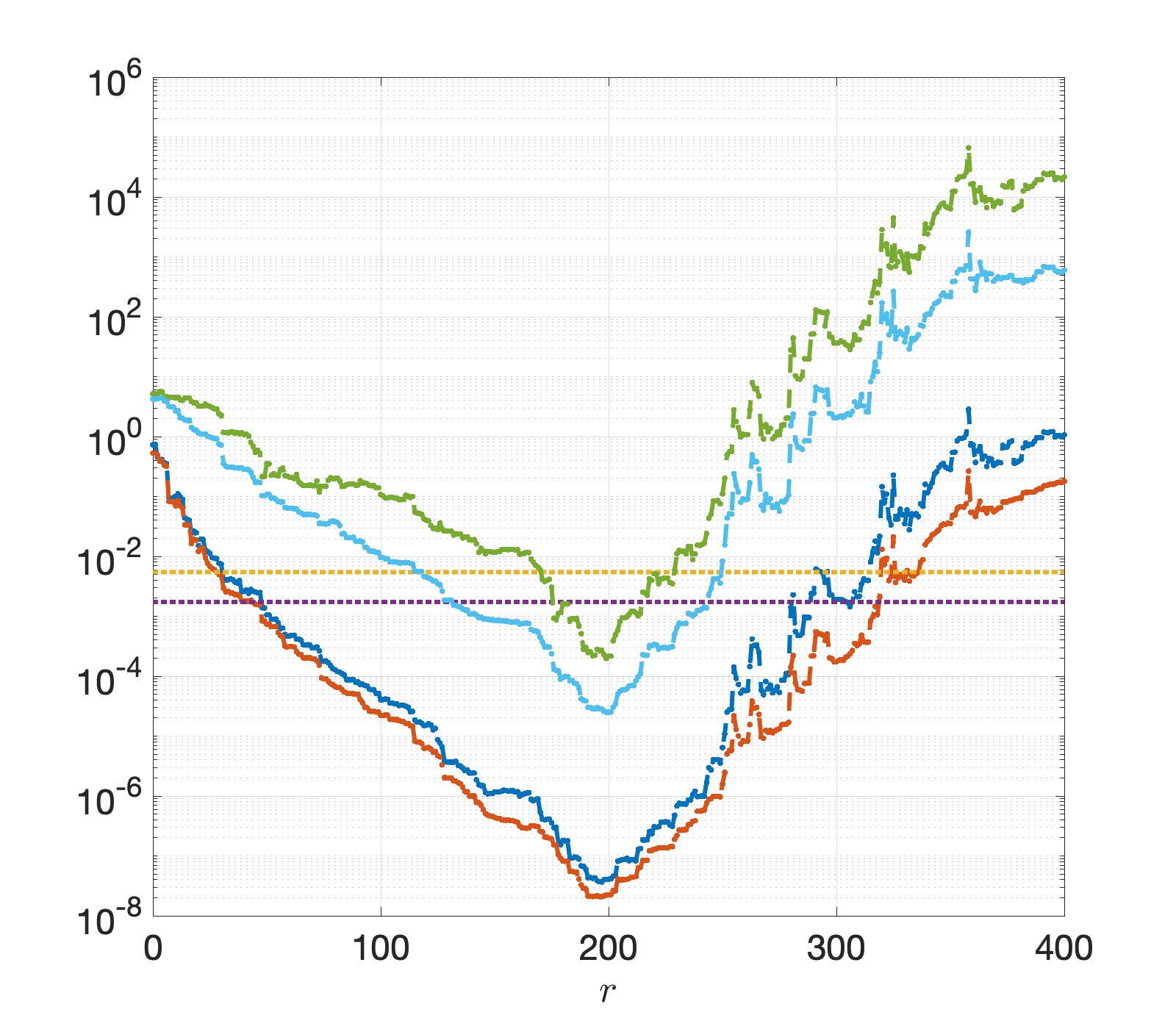}
}
\subfigure[$D = \{2,400,400,1\}$]
{\includegraphics[width=0.31\textwidth]{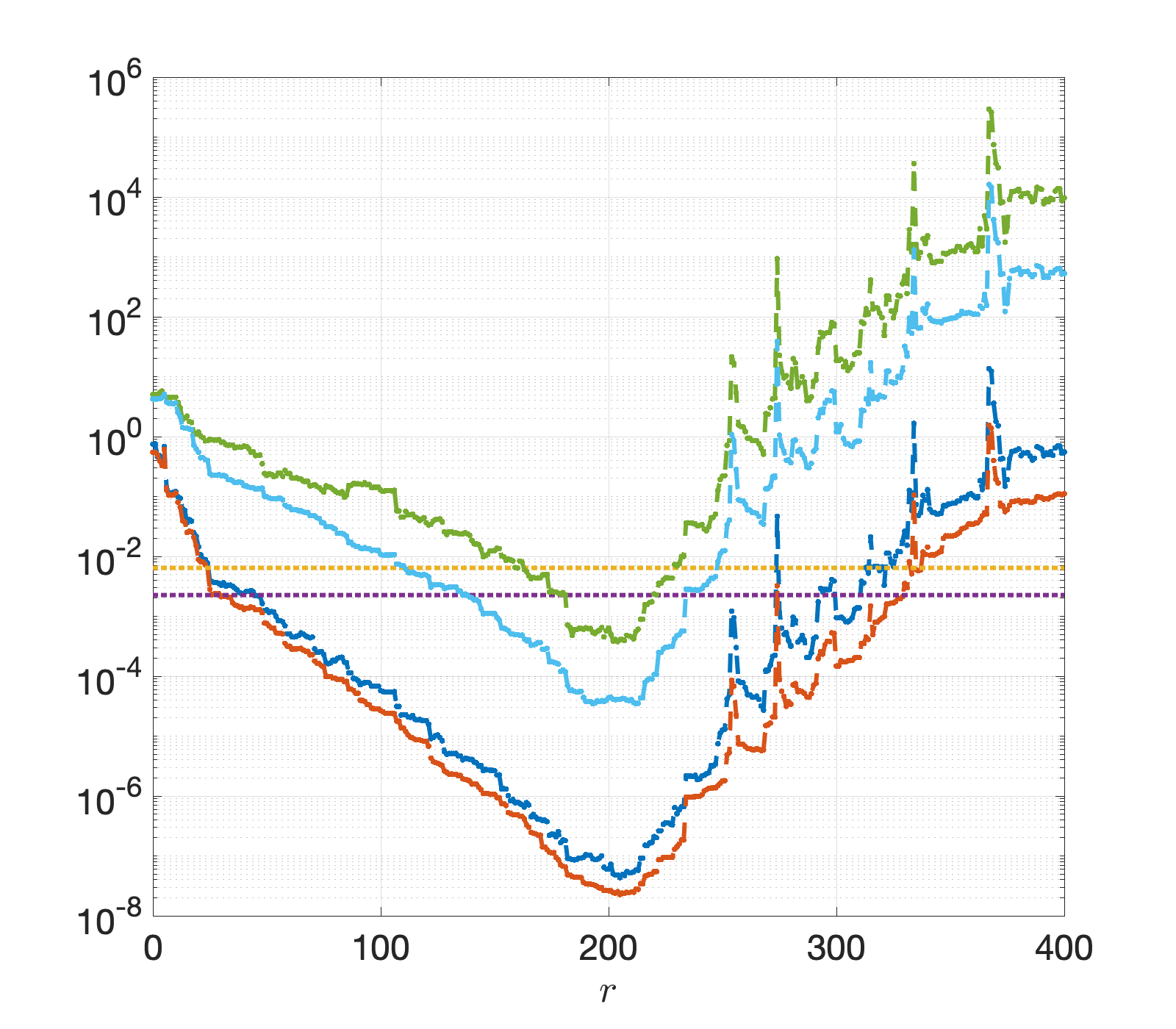}
}
\caption{Poisson equation in the square box $\Omega_2$. The errors $\norm{u_r - u}$ for different values of $r$, measured in $L^{\infty}$ (blue) and $L^2$ (red), for problem (ii), across PINN architectures. We also show the PINN errors $\norm{u^{\theta}_\text{NN} - u}$ in $L^{\infty}$ (yellow) and $L^2$ (purple) as the horizontal curves. The residuals $\norm{e_r}$ in in $L^{\infty}$ (green) and $L^2$ (cyan) are also shown; note the close similarity in the behaviour of the error and residual curves.}\label{PI_Poi2}
\end{figure}

\begin{figure}[tbph]
\centering
\subfigure[$D = \{2,200,400,1\}$]
{\includegraphics[width=0.31\textwidth]{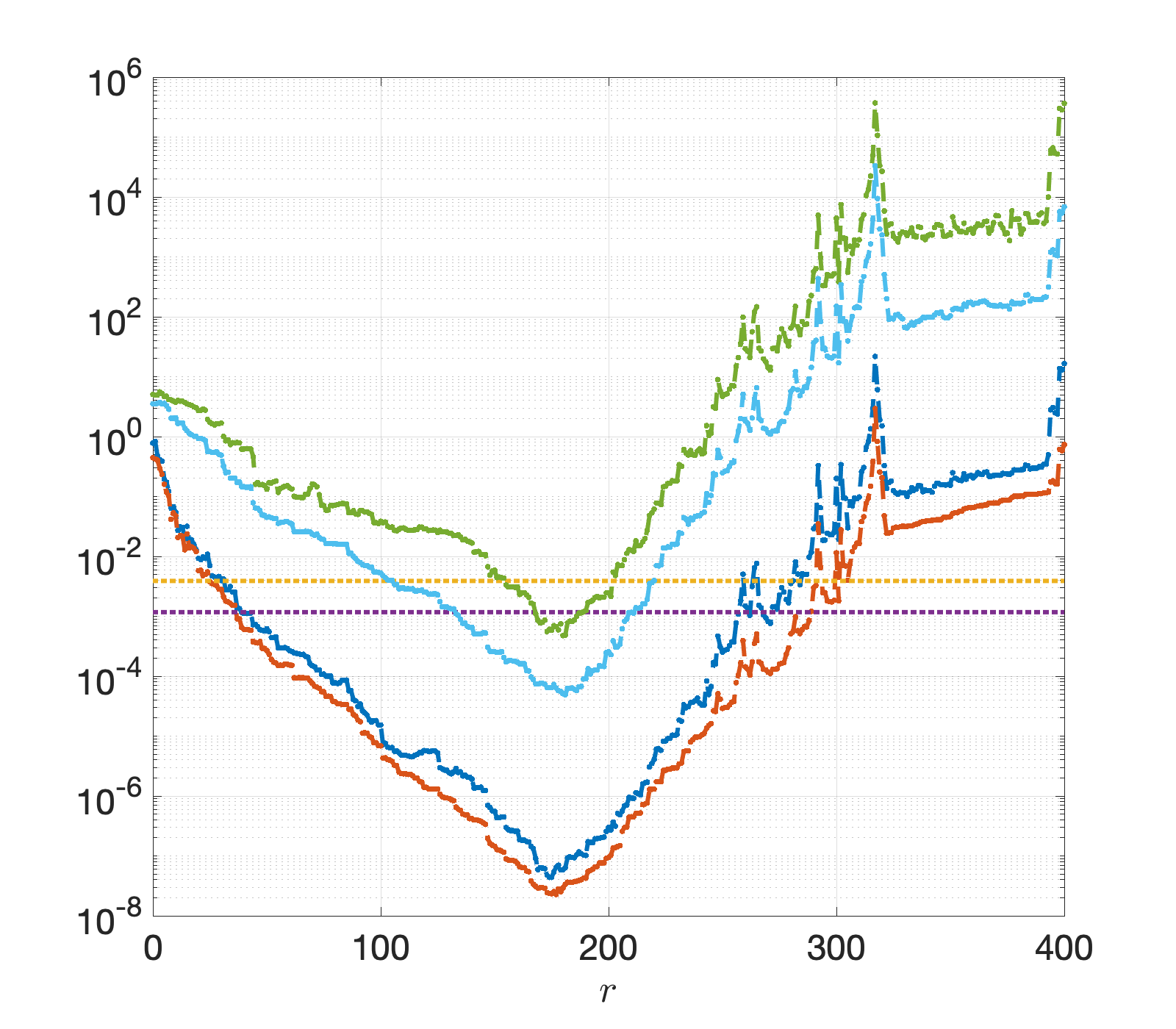}
}
\subfigure[$D = \{2,300,400,1\}$]
{\includegraphics[width=0.31\textwidth]{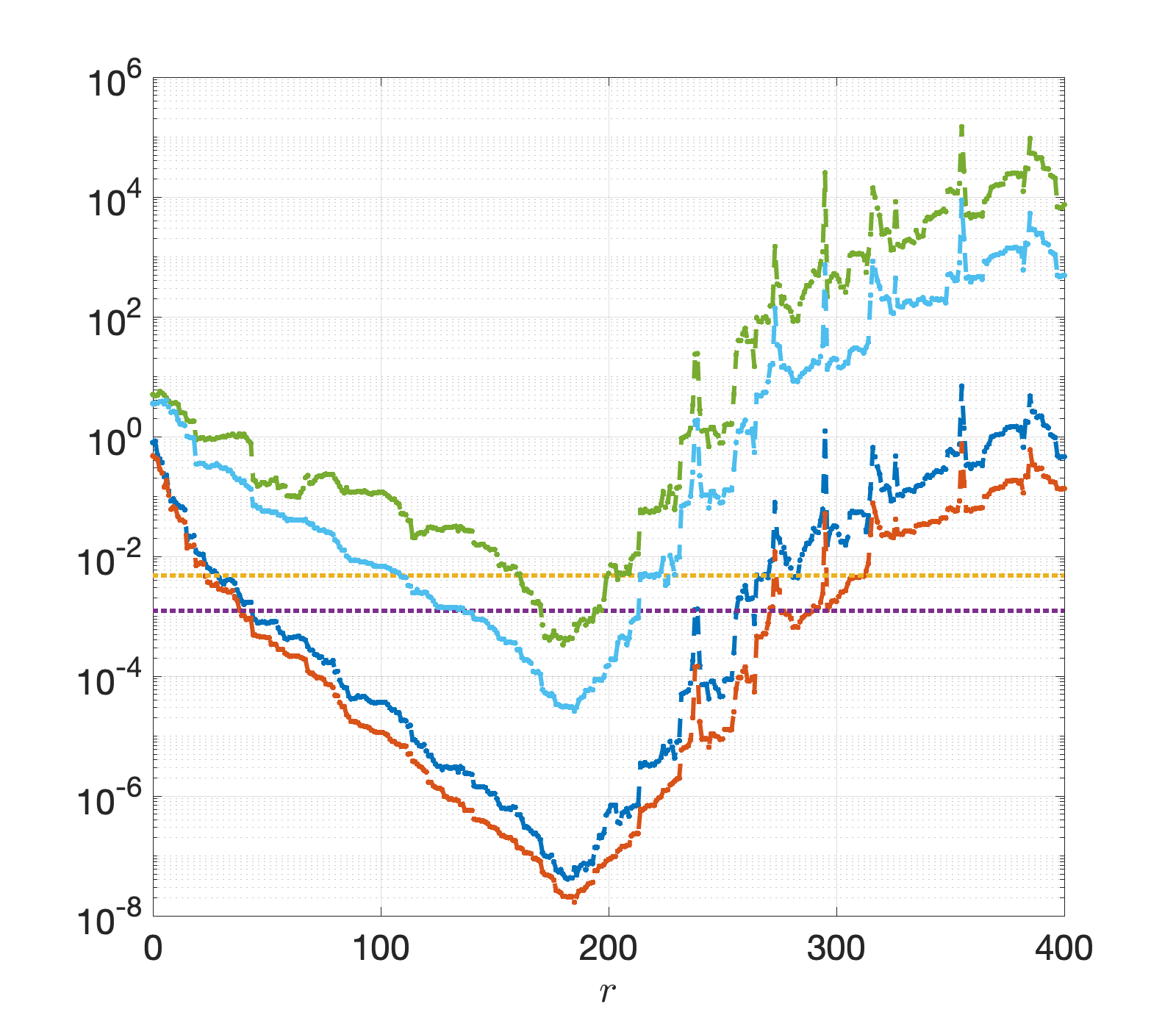}
}
\subfigure[$D = \{2,400,400,1\}$]
{\includegraphics[width=0.31\textwidth]{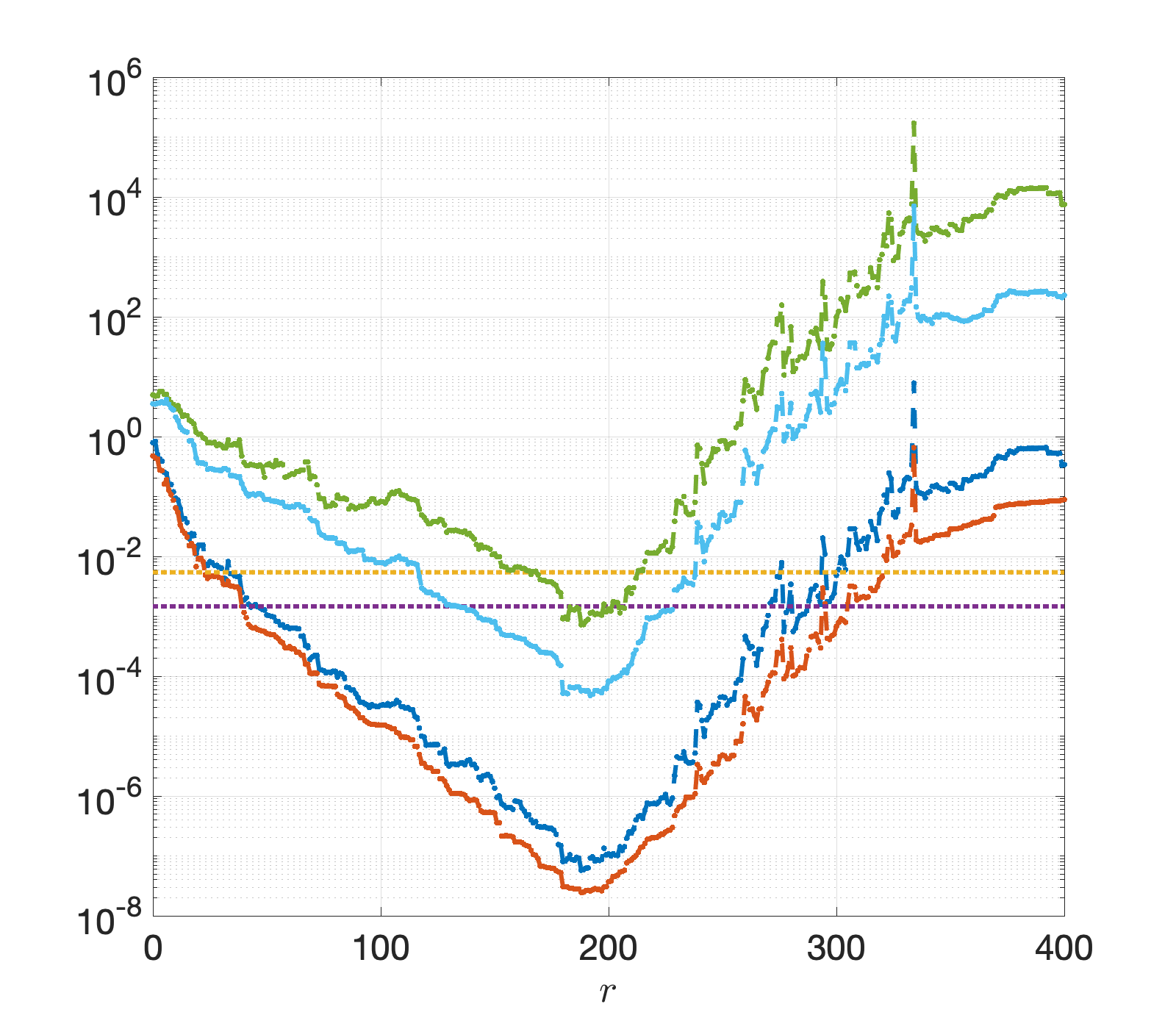}
}
\caption{Poisson equation on the L-shaped domain $\Omega_3$. The errors $\norm{u_r - u}$ for different values of $r$, measured in $L^{\infty}$ (blue) and $L^2$ (red), for problem (iii), across PINN architectures. We also show the PINN errors $\norm{u^{\theta}_\text{NN} - u}$ in $L^{\infty}$ (yellow) and $L^2$ (purple) as the horizontal curves. The residuals $\norm{e_r}$ in in $L^{\infty}$ (green) and $L^2$ (cyan) are also shown; note the close similarity in the behaviour of the error and residual curves.}\label{PI_PoiL}
\end{figure}

We use the Tanh activation function and train the PINNs using Adam with a learning rate of $10^{-3}$ for $5000$ epochs for (i) and 10000 for (ii) and (iii). The linear systems are solved on a high-order Gauss--Legendre quadrature grid, with a finer grid used for calculating the errors. The results are shown in Figures \ref{PI_Poi1}, \ref{PI_Poi2} and \ref{PI_PoiL} for various layer widths $D = \{d_0,d_1,\hdots,d_{L+1}\}$. Note that the profiles of the errors $\norm{u_r - u}$ exhibit a V-shape: they decay rapidly for small to moderate $r$, before hitting a trough, and start rising for larger $r$. As noted earlier, this is a consequence of the the division by small singular values in \eqref{CoB}. Regardless, we note that the lowest errors are typically four to five orders of magnitude smaller than the PINN errors. 

We also note that the residual profiles $\norm{e_r}$ bear close resemblance to the error curves. This serves to establish the use of the former as a useful proxy for setting $r$ optimally: by choosing it so as to minimize $\norm{e_r}$, we invariably obtain the best approximation. We also draw attention to the fact that the optimal values of $r$ rise for larger networks, indicating the increasing expressivity of the corresponding space $\mathcal{S}_{\theta}$.

\section{Extension to time-dependent and nonlinear problems}\label{SecExtend}

The approach described in the previous section can be extended to problems other than the Poisson equation by training PINNs to solve them, calculating bases for the spaces defined by last hidden layers, and using them in variational formulations. Alternatively, we can employ the bases obtained from the Poisson equation over the various domains and use them directly. This has the benefit of not requiring problem-specific training as well as allowing us to explore the transfer learning capabilities of these functions. 

As an example, consider the model equation
\eqn{
u_t(t,x) &=& \nabla \cdot (k(x) \nabla u(t,x)) + \mathcal{N}[u] + f(t,x), \quad x \in \Omega, \ t > 0, \nonumber\\
u(t,x) &=& g(t,x), \quad x \in \partial\Omega,  \ t > 0, \nonumber\\
u(0,x) &=&u_0(x), \quad x \in \Omega, \label{NLTEqn}
}
where $\mathcal{N}$ is some nonlinear operator. We employ the following iteration scheme, obtained from the four-step Backward Differentiation Formula (BDF-4), to discretize the time derivative \cite{hundsdorfer2007imex}
\eqn{
&& u^{n+1} - \frac{12\Delta t}{25} \left(\nabla \cdot (k(x) \nabla u^{n+1})\right) = \frac{1}{25}\left(48u^{n} - 36u^{n-1} + 16u^{n-2} - 3u^{n-3}\right) +\nonumber\\
&& \hspace{4cm} \frac{12\Delta t}{25}\left( f^{n+1}  + 4\mathcal{N}\left[u^n\right] - 6\mathcal{N}\left[u^{n-1}\right] + 4\mathcal{N}\left[u^{n-2}\right] - \mathcal{N}\left[u^{n-3}\right]\right), \label{bdf4}
}
where $u^{n}(x) \approx u(t_n,x)$ and $f^{n+1}(x) = f(t_n,x)$. Coupled with the boundary conditions $u^{n+1}(x) = g((n+1)\Delta t,x)$ for $x \in \partial \Omega$, this problem can be solved iteratively with simple modifications to the linear system \eqref{Abdefn}. Note that the matrix $A$ needs to be built just once (for a specified $\Delta t$) for the entirety of a simulation. The use of a high-order implicit-explicit scheme confers better accuracy and stability properties and allows us to avoid the time-step restrictions due to the stiffness of the diffusion term. However, it does require the values of the solution for the first three steps to be provided, or a one-step routine for jump-starting. In the results shown here, we have used the backward Euler scheme for this purpose. 

\begin{figure}[tbph]
\centering
\subfigure[$D = \{1,30,30,1\}$]
{\includegraphics[width=0.31\textwidth]{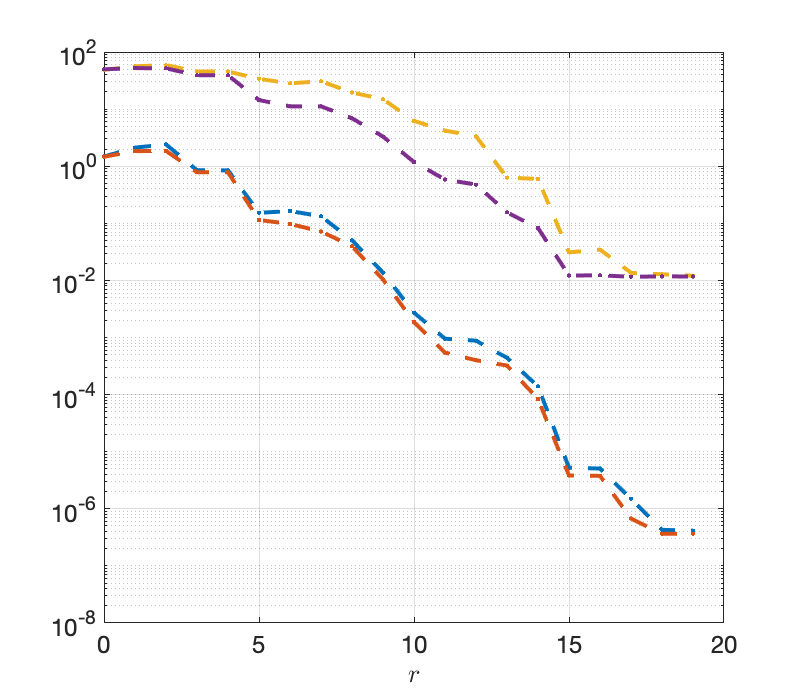}
}
\subfigure[$D = \{1,60,60,1\}$]
{\includegraphics[width=0.31\textwidth]{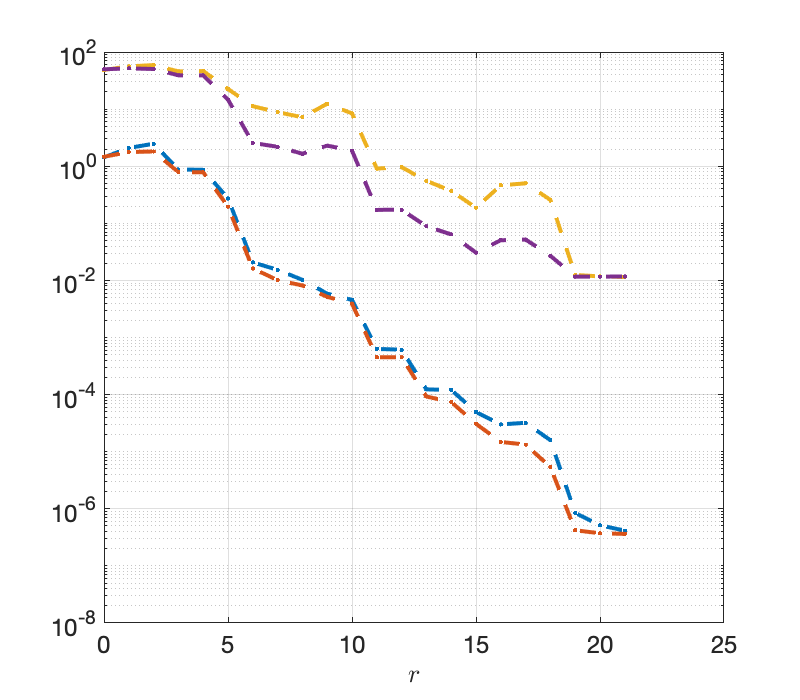}
}
\subfigure[$D = \{1,50,50,50,1\}$]
{\includegraphics[width=0.31\textwidth]{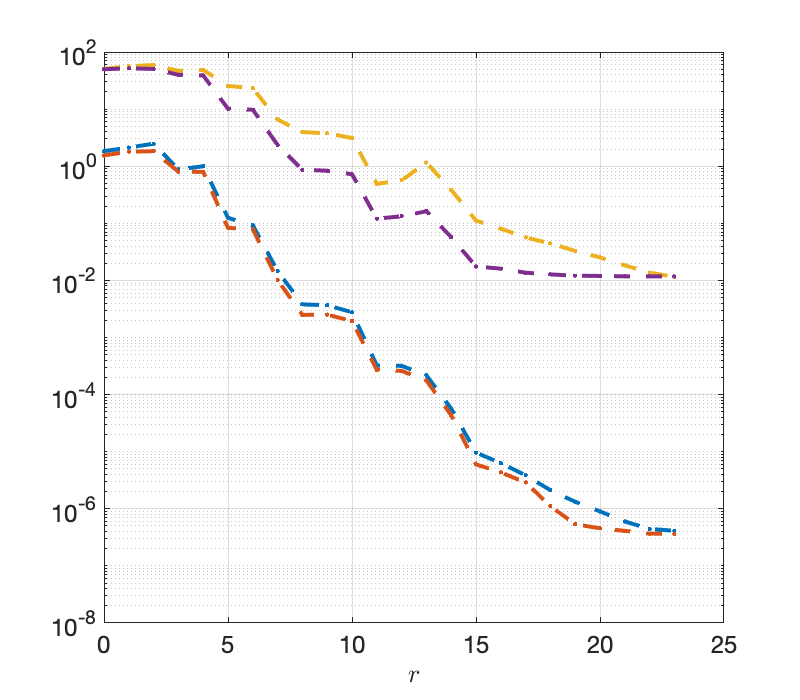}
}
\caption{The heat equation on $\Omega_1 = [-1,1]$; we show the errors $E_{r,p}$ for different values of $r$, measured with $p = \infty$ (blue) and $p = 2$ (red). We also show present the $L^\infty$ (yellow) and $L^2$ (purple) norms of the residuals.}\label{PI_Heat1}
\end{figure}

\begin{figure}[tbph]
\centering
\subfigure[$D = \{2,200,400,1\}$]
{\includegraphics[width=0.31\textwidth]{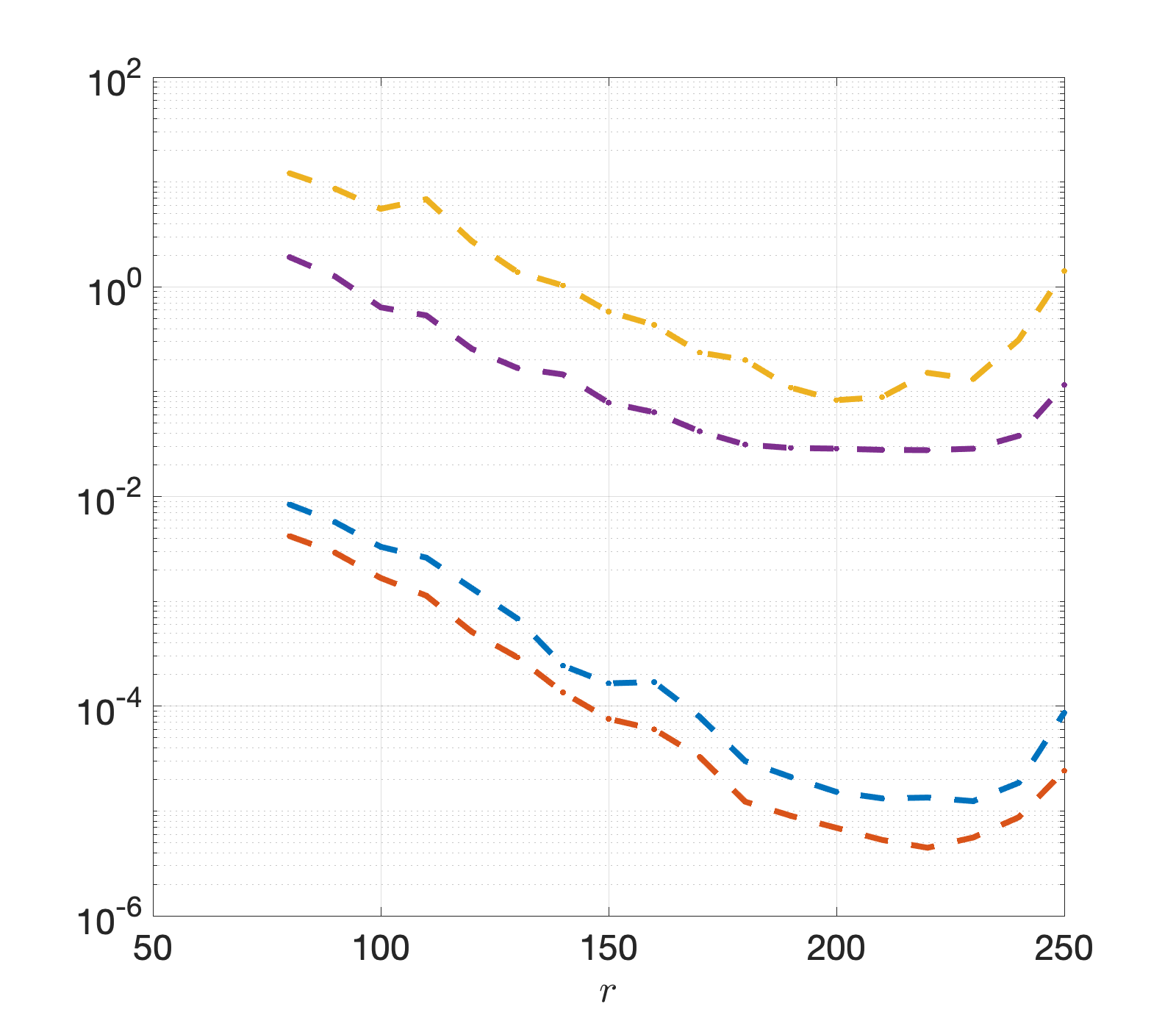}
}
\subfigure[$D = \{2,300,400,1\}$]
{\includegraphics[width=0.31\textwidth]{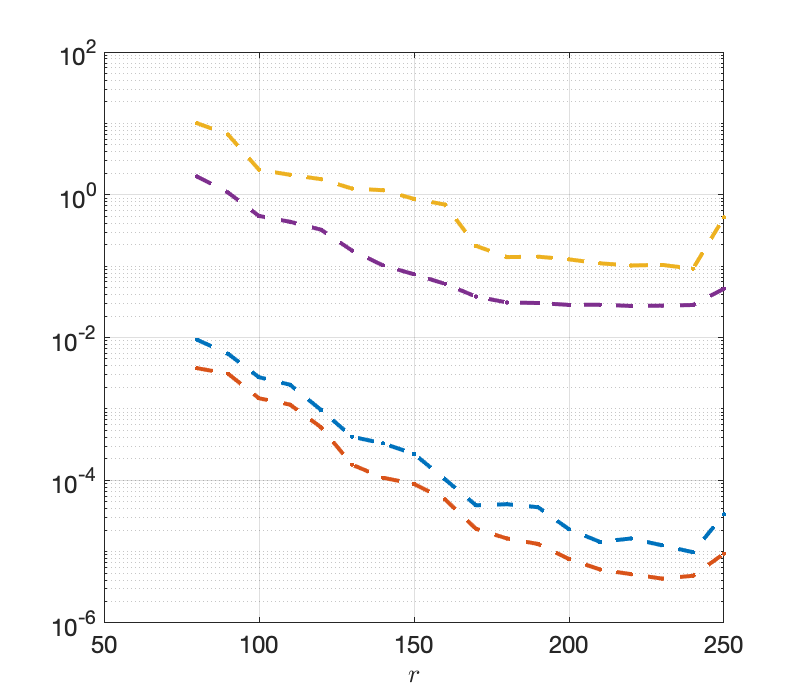}
}
\subfigure[$D = \{2,400,400,1\}$]
{\includegraphics[width=0.31\textwidth]{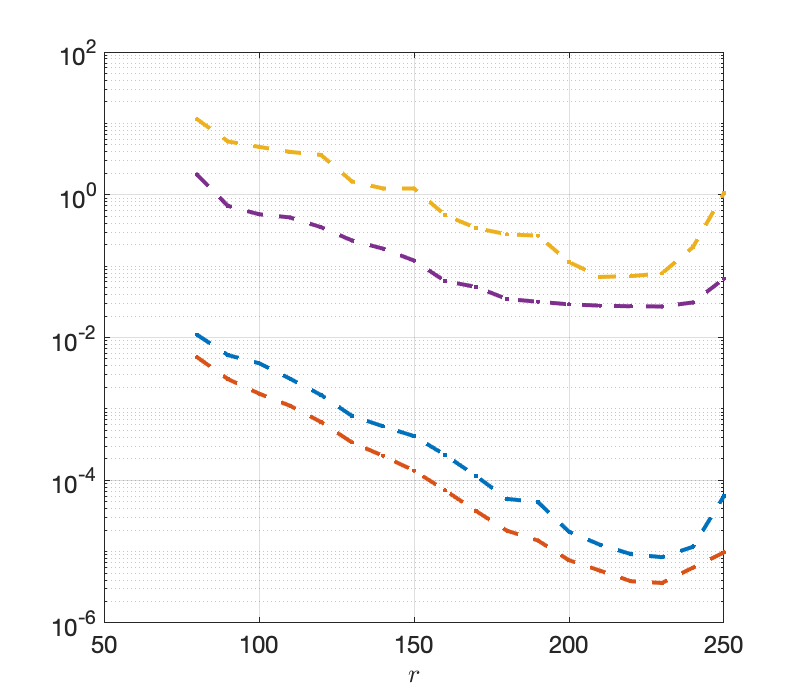}
}
\caption{The heat equation on the square box $\Omega_2$; the errors $E_{r,p}$ are shown for different values of $r$, measured with $p = \infty$ (blue) and $p = 2$ (red), as well as the averaged $L^\infty$ (yellow) and $L^2$ (purple) norms of the residuals.}\label{PI_Heat2}
\end{figure}

\begin{figure}[tbph]
\centering
\subfigure[$D = \{2,200,400,1\}$]
{\includegraphics[width=0.31\textwidth]{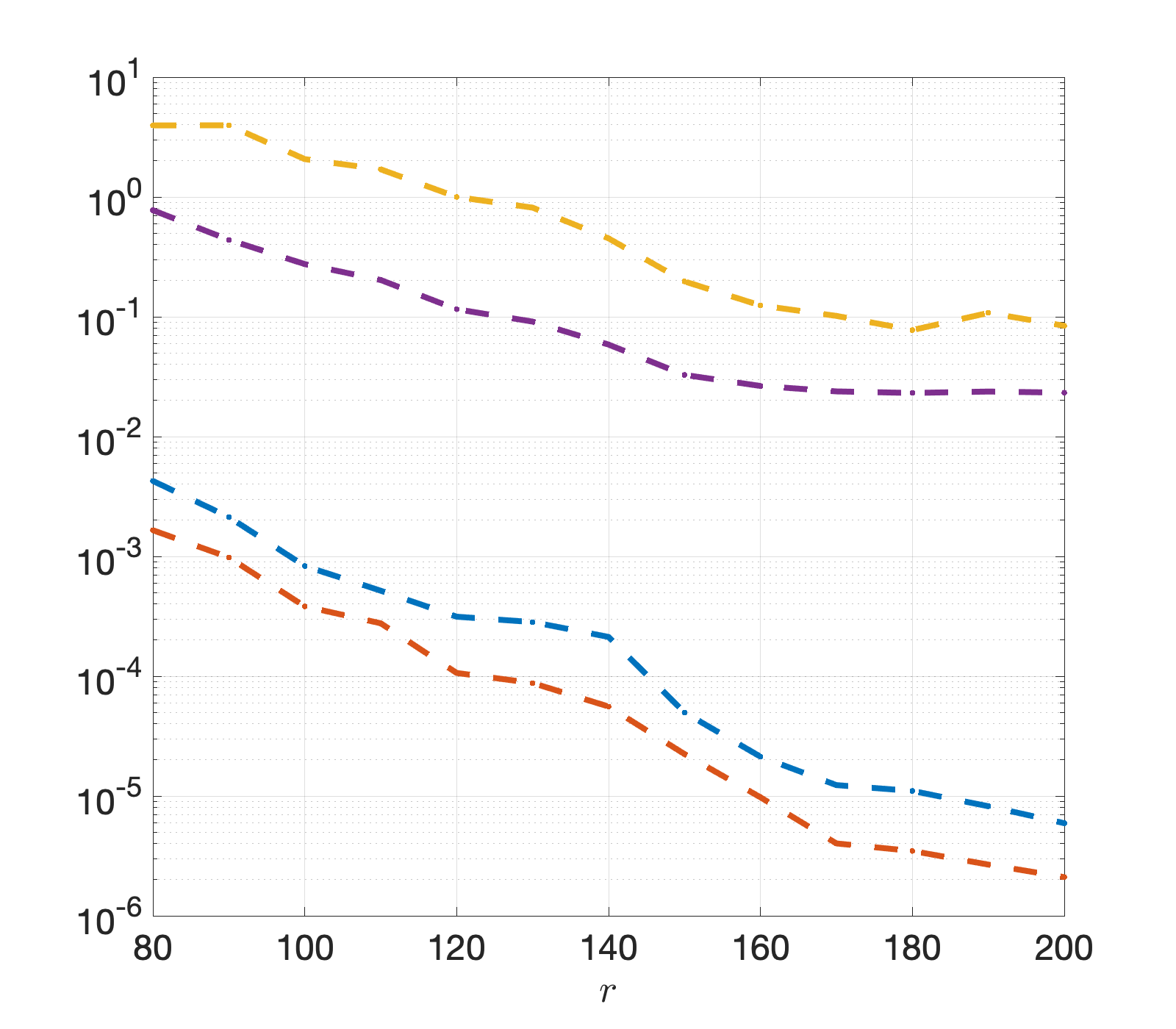}
}
\subfigure[$D = \{2,300,400,1\}$]
{\includegraphics[width=0.31\textwidth]{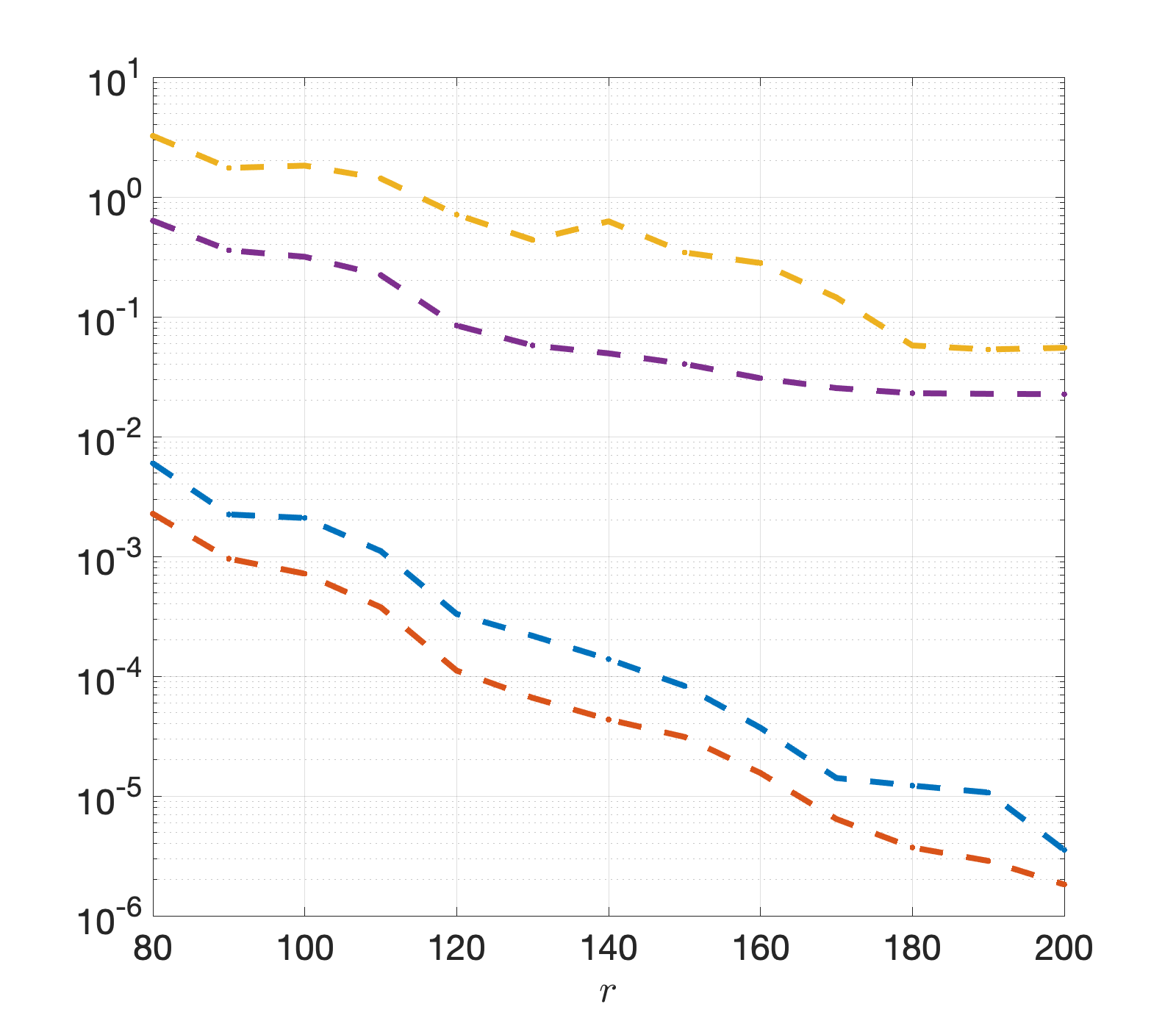}
}
\subfigure[$D = \{2,400,400,1\}$]
{\includegraphics[width=0.31\textwidth]{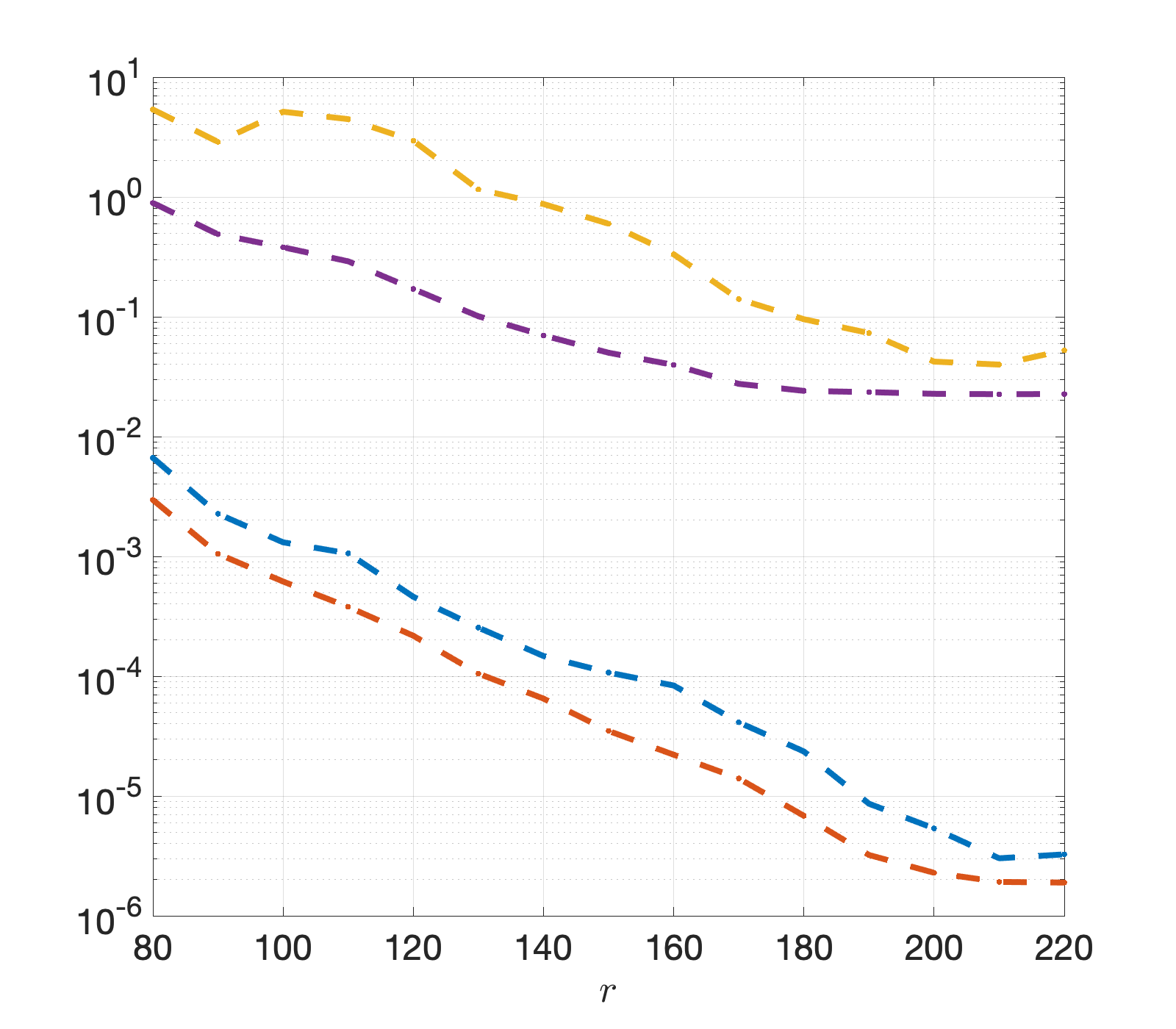}
}
\caption{The heat equation in the L-shaped domain $\Omega_3$. We present the errors $E_{r,p}$ for different values of $r$, measured with $p = \infty$ (blue) and $p = 2$ (red), along with the averaged $L^\infty$ (yellow) and $L^2$ (purple) norms of the residuals.}\label{PI_HeatL}
\end{figure}

To test this approach, we use $\mathcal{N} = 0$ and $k(x) \equiv 1$ (to obtain the heat equation) and the exact solutions $u(t,x) = \sin(5x)\cos(4t)$ for the 1D domain and $u(t,x,y) = \sin(e^x + y^3) e^{\cos(4t)}$ for the two 2D domains to calculate the forcing functions and to compare against the computed solutions to obtain the errors. In Figures \ref{PI_Heat1}, \ref{PI_Heat2}, and \ref{PI_HeatL}, we display the errors
\eqn{
E_{r,p} = \left(\int_0^T \norm{u_r(t,\cdot) - u(t,\cdot)}_{L^p(\Omega)}^2 \ dt \right)^{1/2} \nonumber
}
as they decay with increasing $r$ across architectures, for $p = 2$ and $p = \infty$, using the same basis functions as for the results in Figures \ref{PI_Poi1}, \ref{PI_Poi2}, and \ref{PI_PoiL}. These plots do not show the PINN baselines as the basis functions in use have been computed from a PINN for a different problem.

This framework can also be leveraged to solve stationary nonlinear equations of the form
\eqn{
\nabla \cdot (k(x) \nabla u(x)) + \mathcal{N}[u] + f(x) &=& 0, \quad x \in \Omega, \nonumber\\
u(x) &=& g(x), \quad x \in \partial\Omega. \label{NLEqn}
}

\begin{figure}[tbph]
\centering
\subfigure[$D = \{1,30,30,1\}$]
{\includegraphics[width=0.31\textwidth]{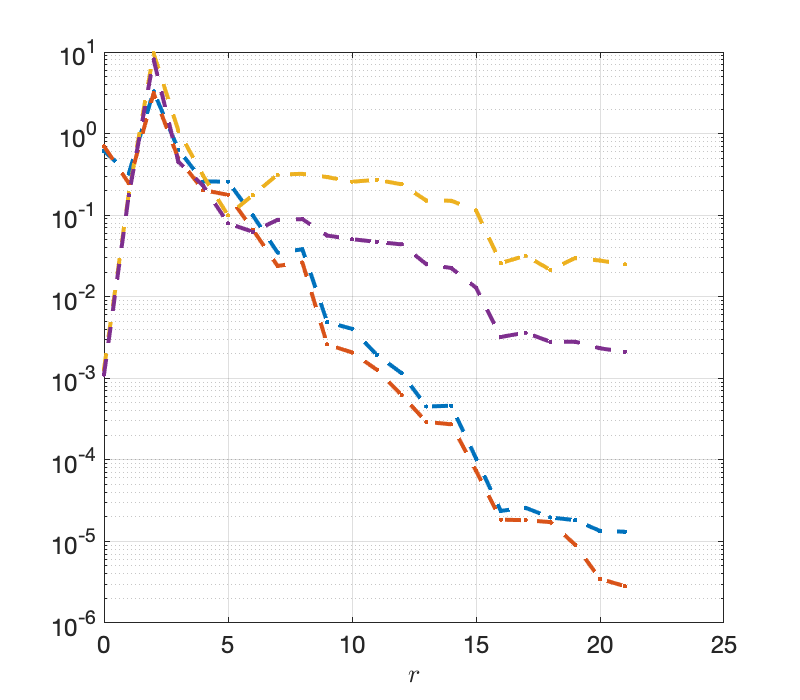}
}
\subfigure[$D = \{1,60,60,1\}$]
{\includegraphics[width=0.31\textwidth]{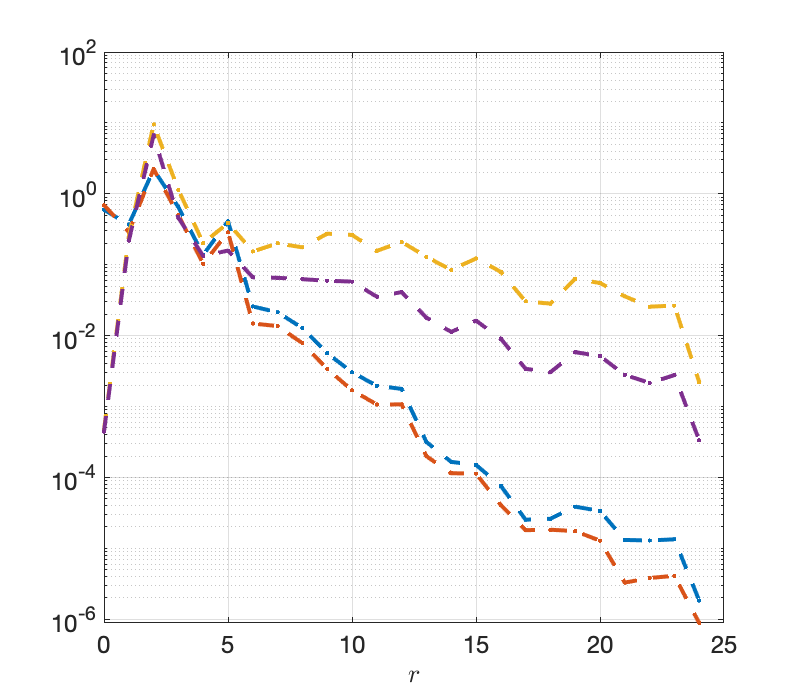}
}
\subfigure[$D = \{1,50,50,50,1\}$]
{\includegraphics[width=0.31\textwidth]{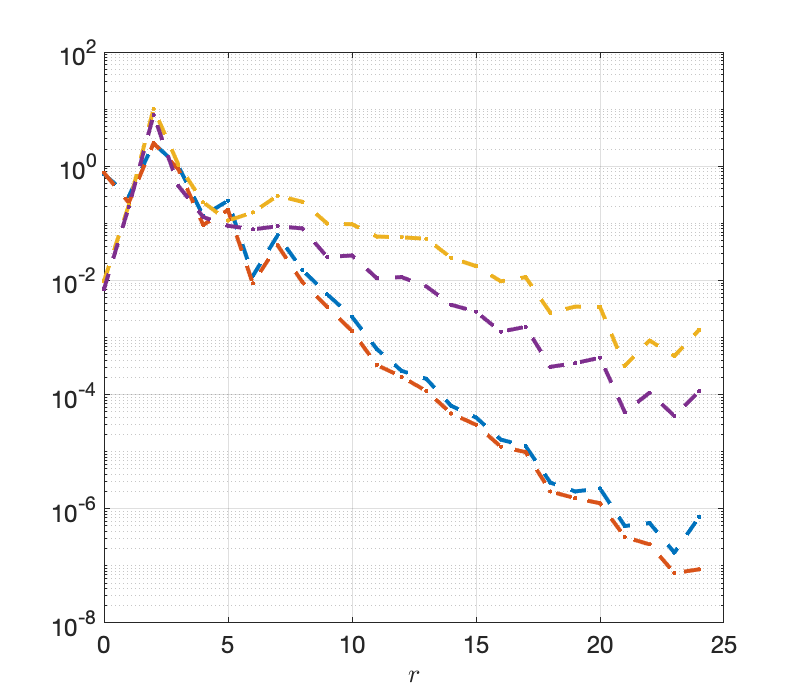}
}
\caption{The steady viscous Burgers' equation \eqref{SVBeqn} on the 1D interval $\Omega_1$. The errors for different values of $r$ are shown, measured in $L^\infty$ (blue) and $L^2$ (red), along with the $L^\infty$ (yellow) and $L^2$ (purple) norms of the residuals.}\label{PI_SVB1}
\end{figure}

By considering \eqref{NLEqn} as the steady-state of \eqref{NLTEqn}, we obtain an iteration scheme that, starting from a suitably chosen initial condition, converges to the solution of \eqref{NLEqn}. In Figure \ref{PI_SVB1}, we display the decay of the norms of the solution errors and residuals with increasing $r$ for the steady viscous Burgers' equation
\eqn{
\nu u''(x) &=& u(x)u'(x),  \quad x \in (-1,1), \nonumber\\
u(\pm1) &=& g_{\pm}, \label{SVBeqn}
}
with the exact profile $u(x) = -0.6\tanh\left(\frac{0.3}{\nu}(x - 0.2)\right)$ used for boundary data and for comparison. We set $\nu = 0.1$ and solve the time-dependent problem till $T = 400$ with $\Delta t = 0.1$, with a linear function that obeys the boundary conditions set as the initial profile. Observe that the errors decay as for the earlier problems, with the residual again serving as a guide for the optimal choice of $r$.

In Figure \ref{PBeqn}, we show the results for the Poisson--Boltzmann equation
\eqn{
-u''(x) + K^2 \sinh(u) &=& 0, \quad x \in (-1,1), \nonumber\\
u(\pm 1) &=& g_{\pm}. \label{PBeqn} 
}

\begin{figure}[tbph]
\centering
\subfigure[$D = \{1,30,30,1\}$]
{\includegraphics[width=0.31\textwidth]{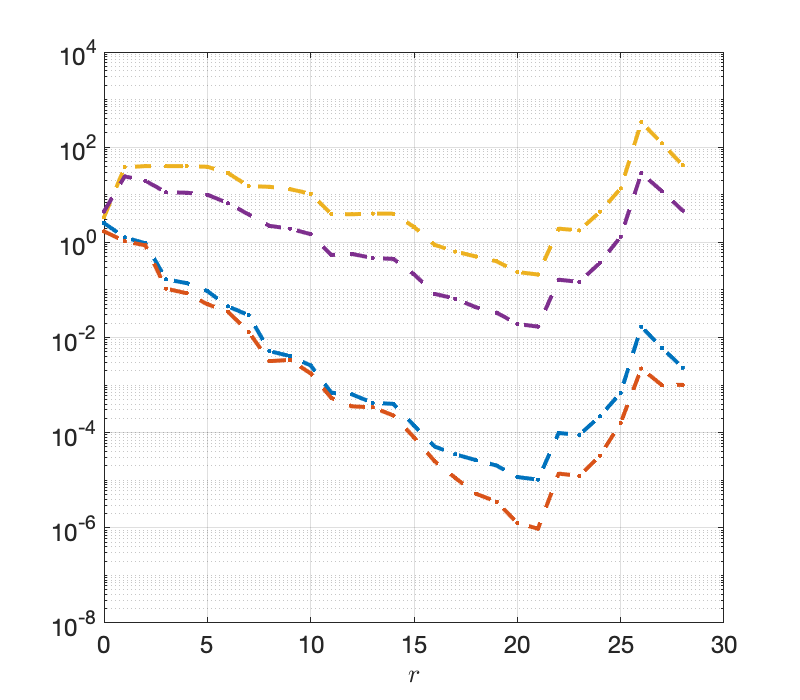}
}
\subfigure[$D = \{1,60,60,1\}$]
{\includegraphics[width=0.31\textwidth]{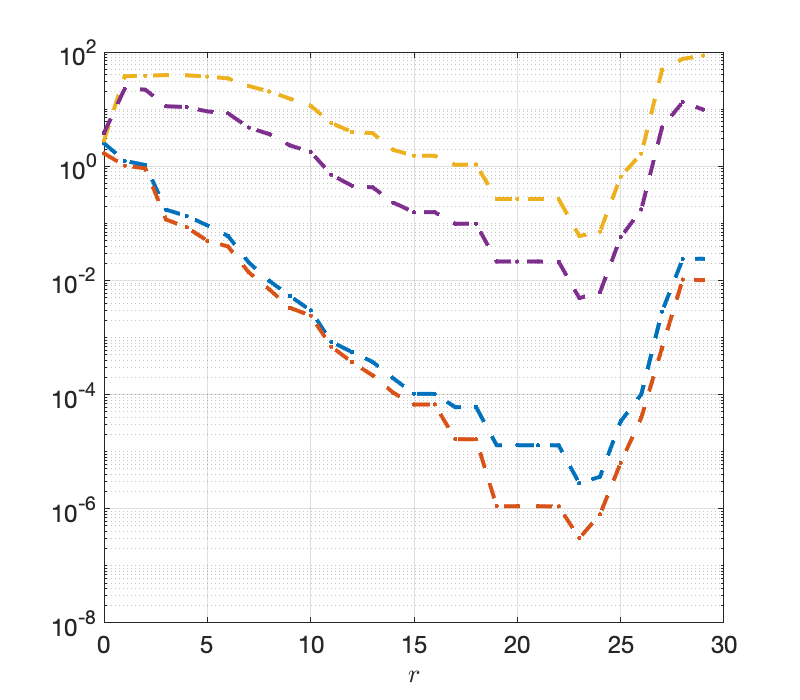}
}
\subfigure[$D = \{1,50,50,50,1\}$]
{\includegraphics[width=0.31\textwidth]{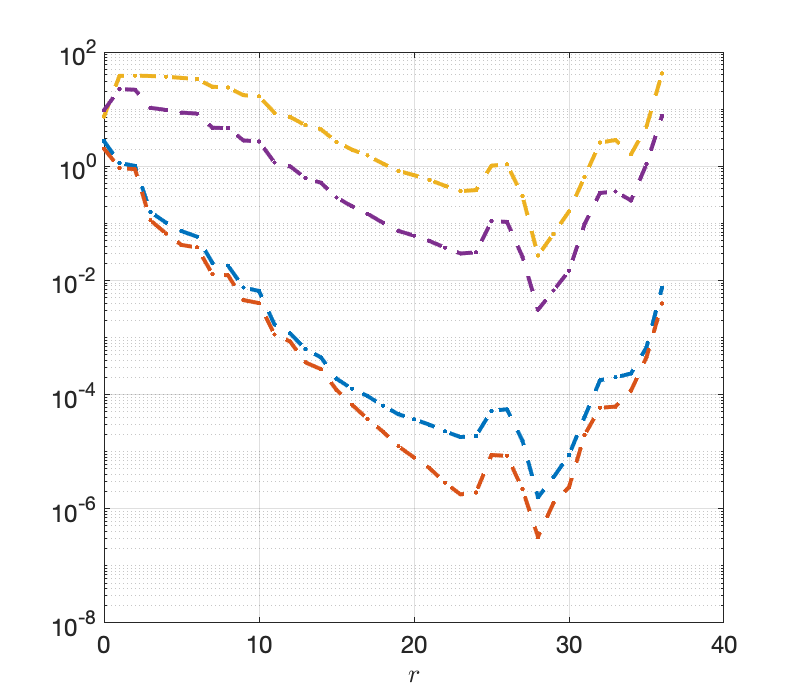}
}
\caption{The Poisson--Boltzmann equation \eqref{PBeqn} on the 1D interval $\Omega_1$. The errors for different values of $r$ are shown, measured in $L^\infty$ (blue) and $L^2$ (red), along with the $L^\infty$ (yellow) and $L^2$ (purple) norms of the residuals.}\label{PI_PB1}
\end{figure}

We compare the computed solutions to those yielded by using Legendre polynomials for $g_+ = 2$ and $g_- = -3$. We display the results for $K = 2$ in Figure \ref{PI_PB1}, solving the time-dependent problem till $T = 100$ with $\Delta t = 0.1$. We again find that the errors and residuals decay fairly similarly, further illustrating the cross-problem adaptability of the basis functions yielded by our approach.

\section{Discussion}\label{SecDisc}

Our preliminary work demonstrates that using a trained PINN to extract basis functions for use in a spectral method can significantly boost the accuracy of the parent network. Due to its simplicity, it can be applied as a post-processing step at the end of training any such architecture. Moreover, the use of \eqref{CoB} for differentiation and for evaluation away from the quadrature grid allows us to bypass polynomial projections as in \cite{meuris2023machine}, while also providing a residual-based metric that can be used to choose the optimal number of basis functions.

The explicit identification of basis functions orthonormal with respect to a high-order quadrature rule makes them well-suited for use in variational formulations. The resulting linear systems are small and well-conditioned, in contrast with least-squares problems yielded by collocation approaches. We also demonstrate that the basis functions can be employed to solve more complex problems on the same domains, whether they are time-dependent and/or possess different nonlinear mechanisms. This attests to the generalizability of the functions obtained from our approach.

\section{Acknowledgments}\label{SecAck}
The work of SQ and PS is supported by the Department of Energy (DOE) Office of Science, Advanced Scientific Computing Research program under Resolution-invariant Deep Learning for Accelerated Propagation of Epistemic and Aleatory Uncertainty in Multi-scale Energy Storage Systems, and Beyond (Project No. 81824). Pacific Northwest National Laboratory is a multi-program national laboratory operated for the U.S. Department of Energy by Battelle Memorial Institute under Contract No. DE-AC05-76RL01830.

\bibliography{refs}

\end{document}